\documentclass[11pt]{article}

\usepackage{amsmath,amsthm,amssymb}
\usepackage{enumerate}
\usepackage{graphicx}
\usepackage{epsfig}

\usepackage{hyperref}
\hypersetup{
  colorlinks=true,
  linkcolor=blue,
  filecolor=magenta,
   urlcolor=blue,
}
\usepackage{enumitem}
\usepackage{hyperref}

 \usepackage{geometry}
\newcommand*{\mailto}[1]{\href{mailto:#1}{\nolinkurl{#1}}}

\newtheorem{theorem}{Theorem}[section]
\newtheorem{definition}[theorem]{Definition}
\newtheorem{lemma}[theorem]{Lemma}
\newtheorem{example}[theorem]{Example}
\newtheorem{proposition}[theorem]{Proposition}
\newtheorem{corollary}[theorem]{Corollary}

\newtheorem{remark}[theorem]{Remark}

\newcommand{\fr}{\frac}

\newcommand{\R}{{\mathbb R}}

\newcommand{\Z}{{\mathbb Z}}
\newcommand{\Co}{{\mathbb C}}

\newcommand{\cA}{{\cal A}}

\newcommand{\st}{\stackrel}
\newcommand{\toEw}{\st{E}{-\!\!-\!\!\!\rightharpoonup}}
\newcommand{\toE}{\st{E}{-\!\!\!-\!\!\!\to}}

\newcommand{\toHs}{\st{H^s}{-\!\!\!-\!\!\!\to}}

\newcommand{\toCbE }{\st{C_b(\R,E)}{-\!\!\!-\!\!\!-\!\!\!-\!\!\!-\!\!\!-\!\!\!-\!\!\!\to}}
\newcommand{\toCE}{\st{C(\R,E)}{-\!\!\!-\!\!\!-\!\!\!-\!\!\!-\!\!\!-\!\!\!\to}}
\newcommand{\toCTE}{\st{C([0,T],E)}{-\!\!\!-\!\!\!-\!\!\!-\!\!\!-\!\!\!-\!\!\!-\!\!\!\to}}

\newcommand{\toCHs}{\st{C([-T,T],H^s)}{-\!\!\!-\!\!\!-\!\!\!-\!\!\!-\!\!\!-\!\!\!-\!\!\!-\!\!\!-\!\!\!\to}}
\newcommand{\toCnHs}{\st{C([0,T],H^s)}{-\!\!\!-\!\!\!-\!\!\!-\!\!\!-\!\!\!-\!\!\!-\!\!\!-\!\!\!-\!\!\!\to}}

\newcommand{\toCLr}{\st{C([0,T],L^r)}{-\!\!\!-\!\!\!-\!\!\!-\!\!\!-\!\!\!-\!\!\!-\!\!\!-\!\!\!\to}}
\newcommand{\toC}{\st{C(\R)}{-\!\!\!-\!\!\!-\!\!-\!\!\!\to}}
\newcommand{\cH}{{\cal H}}

\newcommand{\cU}{{\cal U}}

\newcommand{\al}{\alpha}
\newcommand{\om}{\omega}
\newcommand{\vp}{\varphi}

\newcommand{\De}{\Delta}
\newcommand{\de}{\delta}
\newcommand{\vka}{\varkappa}

\newcommand{\ga}{\gamma}
\newcommand{\ve}{\varepsilon}

\newcommand{\Om}{\Omega}

\newcommand{\ti}{\tilde}
\newcommand{\na}{\nabla}
\newcommand{\pa}{\partial}
\newcommand{\dist}{{\rm dist\5}}

\newcommand{\loc}{{\rm loc}}

\newcommand{\rRe}{{\rm Re\5}}
\newcommand{\rIm}{{\rm Im\5}}

\newcommand{\ov}{\overline}

\newcommand{\5}{{\hspace{0.5mm}}}
\newcommand{\ds}{\displaystyle}

\numberwithin{equation}{section}

\newcommand{\ci}{\cite}
\newcommand{\la}{\label}

\newcommand{\be}{\begin{equation}}
 \newcommand{\ee}{\end{equation}}
 \newcommand{\beqn}{\begin{eqnarray}}
 \newcommand{\eeqn}{\end{eqnarray}}
\newcommand{\ba}{\begin{array}}
 \newcommand{\ea}{\end{array}}
\newcommand{\bd}{\begin{definition}}
 \newcommand{\ed}{\end{definition}}
\newcommand{\bt}{\begin{theorem}}
 \newcommand{\et}{\end{theorem}}
\newcommand{\bp}{\begin{proposition}}
 \newcommand{\ep}{\end{proposition}}
\newcommand{\bl}{\begin{lemma}}
 \newcommand{\el}{\end{lemma}}
\newcommand{\bc}{\begin{corollary}}
 \newcommand{\ec}{\end{corollary}}
\newcommand{\bex}{\begin{example}}
 \newcommand{\eex}{\end{example}}
 \newcommand{\br}{\begin{remark} }
 \newcommand{\er}{\end{remark}}
\newcommand{\bce}{\begin{center}}
\newcommand{\ece}{\end{center}}
 \newcommand{\bpr}{\begin{proof}}
\newcommand{\epr}{\end{proof}}

\date{}

\numberwithin{equation}{section}

\begin{document}

\bce
{\huge\bf  On  global attractors  for 2D damped driven
\bigskip\\
 nonlinear Schr\"odinger equations}
\bigskip \bigskip

 {\Large A. I. Komech\footnote{The research supported by the Austrian Science Fund (FWF) under Grant No. P28152-N35.}
 and E. A. Kopylova}
  \medskip
 \\
{\it
\centerline{Faculty of Mathematics of Vienna University
 }
 
{\it
\centerline{Institute for Transmission Information Problems of RAS, Moscow, Russia}
 }

}
 \ece
 \centerline{alexander.komech@univie.ac.at,\quad elena.kopylova@univie.ac.at}

\bigskip\bigskip

\begin{abstract}
Well-posedness and global attractor  are established for 2D damped driven nonlinear Schr\"odinger
equation with almost periodic pumping in a bounded region.
The key role is played by a novel application of the energy equation.

\end{abstract}


\tableofcontents

\section{Introduction}\la{sLM}
We consider the weakly damped driven nonlinear Schr\"odinger
equation in a bounded region $\Om\subset \R^2$ with a smooth boundary $\pa\Om$,
\be\la{dSom}
i\dot\psi(x,t)=-\De\psi(x,t)+f(\psi(x,t))-i\ga\psi(x,t)+p(x,\om t),\qquad x\in\Om,\quad t\in\R,
\ee
where the `friction coefficient' $\ga>0$, and the `pumping' $p(\cdot,\om t)$ is an almost periodic function.
All the derivatives here and below are understood in the sense of distributions.

Equations of type (\ref{dSom}) arise for instance in plasma physics \ci{NB1986} and in optical fibers models \ci{BD1984}.
\smallskip

We impose the Dirichlet boundary condition
\be\la{Dbc}
\psi(x,t)=0,\qquad x\in\pa\Om,\quad t\in\R.
\ee
Our analysis of the nonautonomous Schr\"odinger equation
is motivated by mathematical problems of laser and maser coherent radiation.
In particular,
the maser action is described by the damped driven nonlinear
Maxwell--Schr\"odinger equations in a bounded cavity,
and the Dirichlet boundary condition (\ref{Dbc})
 holds with a high precision  in the case of a solid-state  gain medium \ci{H1984}.

We will identify complex numbers $z\in\Co$ with two-dimensional real vectors $z=(\rRe z,\rIm z)$.
We assume that the nonlinear term is potential, i.e.,
\be\la{fU}
f(\psi)=U'(\psi),\qquad \psi\in\Co\equiv\R^2.
\ee
Here
$U(\psi)\in C^3(\R^2)$ is a real function, and the derivatives mean the derivatives with respect to $\psi_1=\rRe\psi$ and $\psi_2=\rIm\psi$, i.e.,
\be\la{fU2}
f(\psi):= (\pa_1U(\psi), \pa_2U(\psi)),\qquad\psi\in\R^2.
\ee
Respectively, equation (\ref{dSom}) is the identity of $\R^2$-valued distributions,
and the multiplication by $i$ in (\ref{dSom}) means the application of the real
$2\times 2$-matrix
$\left(\ba{rr}0&-1\\1&0\ea\right)$.
Similarly, $f'(\psi)$ and $f''(\psi)$ are tensors:
\be\la{fpr}
f'_{ij}(\psi):=\pa_jf_i(\psi_1,\psi_2)=\pa_j\pa_i U(\psi_1,\psi_2), \qquad f''_{ijk}(\psi):=\pa_j\pa_k f_i(\psi_1,\psi_2),\qquad i,j,k=1,2.
\ee
Note that $f'(\psi)$ is a real symmetric matrix.
The Schr\"odinger equation (\ref{dSom}) under the boundary conditions (\ref{Dbc}) can be written in the `Hamiltonian form'
\be\la{dSHf}
i\dot\psi(t)=D\cH(\psi(t))-i\ga\psi(t)+p(\om t),\qquad t\in\R,
\ee
with the Hamilton functional 
\be\la{Ham}
\cH(\psi):=\int_\Om \fr12|\na\psi(x)|^2\,dx+\cU(\psi),\qquad
\cU(\psi):=\int_\Om U(\psi(x))dx.
\ee
{\bf Our main results.}
For the damped driven Schr\"odinger equation (\ref{dSom}) with almost-periodic pumping in a bounded region
we establish
\medskip\\
I. The well-posedness of the Cauchy problem for equation (\ref{dSom}) (Theorem \ref{pwp});
\\
II. Global attraction in the Sobolev norm $H^1$
of all finite energy solutions to
a compact subset in $H^1$ (Theorem \ref{tA});
\medskip

Let us comment on previous works in these directions. In the case $p(x,t)\equiv 0$ (and $\ga=0$)
the well-posedness for nonlinear Schr\"odinger equations of type (\ref{dS}) was established in \ci{C2003}.
The approach \ci{C2003} relies essentially on the conservation of charge and energy.
\smallskip

Compact attractors in the energy space $H^1$ were constructed
 for weakly damped nonlinear Schr\"odin\-ger equation i) on the circle $\Om=\R/\Z$ or a bounded interval $\Om\subset \R$ by Ghidaglia \ci{G1988},
 ii) on a bounded region $\Om\subset\R^2$ by Abounouh and Goubet \ci{A1993, AG2000}
and iii) on the entire space $\Om=\R^N$ with $N\le 3$ by Lauren\c{c}ot \ci{L1995}.
Similar results were obtained by Ghidaglia and Heron \ci{GH1987} for the Ginzburg--Landau equations on a bounded region $\Om\subset\R^n$ with $n=1,2$.
The methods of these papers rely on the Ball ideas \ci{B2004}.
Note that the pumping terms (or `external force') in \ci{A1993, AG2000, GH1987, L1995} do not depend on time.
In \ci{G1988}, bounded absorbing sets are established
for the
pumping term depending on time,
 while the attractor is constructed for the autonomous case and for time-periodic pumping.
The papers \ci{G1996,G1998,GM2009} are concerned with smoothness properties of functions from
attractors of 1D nonlinear Schr\"odinger equations.

These methods and results were extended i) to 1D {\it nonautonomous}
Schr\"odinger equations on the circle by Wang \ci{W1995}
 and ii) to {\it nonautonomous} KdV and 2D Navier-Stokes equations \ci{MRW2004}.

Tao \ci{T2008} established 
the existence of a~global attractor for radial solutions to
nonlinear autonomous
defocusing Schr\"odinger equation without damping in $\R^n$ with $n\ge 11$.
The paper \ci{S2015} concerns the well-posendess and decay of solutions to 2D damped
{\it autonomous} Schr\"odinger equation in a bounded region.

Recently Cazenave and Han established the long-time bechavior
for nonlinear Schr\"odinger equation in $\R^n$ with a nonlinear subcritical dissipation \ci{CH2020}.
\smallskip

Let us comment on our approach.

First, we prove the existence of a weak solution $\psi(x,t)$, using the standard 
Galerkin approach. Our choice of the dimension two is caused by
the Trudinger inequality, which provides the uniqueness of the weak solution by Theorem 3.6.1 of \ci{C2003}.

Further, we prove that $\psi(x,t)$ is indeed a~strong solution.
This is the first result for {\it nonautonomous}
2D Schr\"odinger equation in a bounded region.
 Note that in our case when $p(x,t)\not\equiv 0$ and $\ga> 0$, the charge and energy conservation do not hold for the nonautonomous equation (\ref{dSom}).
Respectively, the approach of \ci{C2003} is not applicable here.

This is why we introduce a {\it novel method} based on the energy equation \eqref{ebal2}.
This method
allows us to substitute the role of charge and energy conservations from \ci{C2003}
in proving that i) $\psi(x,t)$ is the strong solution and
ii) the solution continuously depends on the initial state and the pumping.

  To prove
the convergence in $H^1$-norm of all finite energy solutions to a compact attractor,
we develop the Ball ideas introduced in the context of autonomous equations \ci{B2004}, and their extensions to nonautonomous 1D Schr\"odinger and KdV equations \ci{W1995,MRW2004}.
The almost-periodicity of the pumping plays a~crucial role in our approach.

\section{A priori estimates for smooth solutions}\la{s2}
Here we prove a priori estimates for sufficiently smooth solutions to the Schr\"odinger equation (\ref{dSom}).
The estimates provide the existence of bounded absorbing sets in the energy norm.
We will get rid of the smoothness assumption in the next section.

To simplify the notation,
we consider equation (\ref{dSom}) with $\om=1$
everywhere in Sections \ref{s2}--\ref{s4} without loss of generality, i.e.,
\be\la{dS}
\dot\psi(x,t)=i\De\psi(x,t)-if(\psi(x,t))-\ga\psi(x,t)-ip(x, t),\qquad x\in\Om,\quad t\in\R.
\ee
Our assumptions are as follows.
\smallskip\\
{\bf Nonlinearity.}
We will assume that the potential $U(\psi)\in C^3(\Co)$, satisfies the following estimates
\begin{gather}
  \vka_1|\psi|^4-b_1\le U(\psi)\le C|\psi|^4,\la{U1}
\\
\rRe f(\psi)\ov \psi\ge\vka_2 U(\psi)-b_2,\la{U2}
\\
   |f(\psi)|
\le \vka_3(1+|\psi|^3), \la{U3}
\\
  | f'(\psi)|\le \vka_4(1+|\psi|^2),\la{U4}
\\
  f'(\psi)\ge -b_3,\la{U5}
\\
  |f''(\psi)|\le \vka_5(1+|\psi|),\la{U6}
\end{gather}
with all $\vka_l>0$ and all $b_l\in\R$.
Let us recall that $f'(\psi)$ in \eqref{U4} and \eqref{U5} denotes the matrix (\ref{fpr}).
For example, the potentials
\be\la{Up}
U(\psi)=a_2|\psi|^4+a_1|\psi|^2+a_0
\ee
satisfy all these conditions if $a_k\in\R$ and $a_2>0$.
Any smooth potential which differs from (\ref{Up}) only in a bounded domain $|\psi|\le C$ also satisfies these conditions.
\smallskip\\
Denote $H^s:=H^s(\Om)$, $L^p:=L^p(\Om)$. Let us introduce the Hilbert phase space
\be\la{E0}
E:= H^1_0(\Om),\qquad \Vert \psi\Vert_{E}=\Vert\na\psi\Vert,
\ee
where $\Vert\cdot\Vert$ stands for the norm in $L^2$.
Then $E^*=H^{-1}$.
\smallskip\\
{\bf The pumping.}
We will assume that
\be\la{pT}
p\in C_b(\R,E)\,\,\,\mbox{\it is a uniform almost periodic function of $t\in\R$,}\quad
p_0:=\Vert p\Vert_{C_b(\R,E)}<\infty.
\ee
By definition (see Appendix of \ci{H1980} for the case of scalar functions), a function $p\in C_b(\R,E)$ is  {\it uniformly almost periodic} if for any sequence $t_k\in\R$ there is a subsequence $t_{k^*}$ such that
the translations $p_{k^*}(t):=p(t+t_{k^*})$ are uniformly  converging,
\be\la{alper}
p_{k^*} \toCbE p^*,\qquad k^*\to\infty.
\ee

\subsection{Energy estimates}
First we prove a priori estimates in the norm of the energy space $E$.
\bl\la{l1}
Let conditions {\rm \eqref{U1}--\eqref{U5}} and {\rm\eqref{pT}} hold,
and $\psi(x,t)\in C^2(\ov\Om\times\ov\R)$ be a solution to \eqref{dS} with $\psi(\cdot,0)\in E$. Then
\be\la{Eb5}
\Vert\psi(t) \Vert_E^2 \le C_0 e^{-\al_\pm t}+D_0,\qquad \pm t \ge 0,
\ee
where $C_0=C_0(\Vert\psi(0)\Vert_E)$, $\al_\pm>0$, and $D_0$ depends on $p_0$ only.
\el
\begin{proof}
Differentiating the Hamilton functional $\cH(t):=\cH(\psi(t))$
and using (\ref{dSHf}), we obtain
the energy balance
\be\la{ebal}
\dot\cH (t)=\langle D\cH(\psi(t)),\dot\psi(t)\rangle=\langle D\cH(\psi(t)),-i D\cH(\psi(t))\rangle+\langle - \De \psi+f(\psi),-\ga \psi(t)-i p(t)\rangle,\qquad t\in\R,
\ee
where the brackets mean the duality, which is the inner product in the {\it real Hilbert space}
$L^2(\Om)\otimes\R^2$.
Hence the first term on the right-hand side vanishes. Therefore,
\begin{align}
\dot\cH (t)&=\langle - \De \psi(t)+f(\psi(t)),-\ga \psi(t)-i p( t)\rangle
\nonumber\\
&=-\ga\langle \na \psi(t),\na\psi(t)\rangle-\ga\langle f(\psi(t)), \psi(t)\rangle
-\langle \na \psi(t),i\na p(t)\rangle-\langle f(\psi(t)), ip(t)\rangle,\quad t\in\R.
\end{align}\la{Eb2}
According to (\ref{U2}),
\be\la{Eb3}
\langle f(\psi(t)), \psi(t)\rangle\ge \vka_2\,{\cal U}(t)-b_2|\Omega|,\quad {\cal U}(t):=\int_{\Omega} U(\psi(x,t))dx.
\ee
Further, \eqref{U1}, (\ref{U3}), and \eqref{pT} imply that
\beqn \nonumber\la{fp3}
\!\!\!\!\!\!\!\!\!\!
|\langle f(\psi(t)), ip(t)\rangle| \!\!&\!\!\le\!\!&\!\!\vka_3\int (1+|\psi(x,t)|^{3})|p(x, t)| dx\\
\!\!&\!\!\le\!\!&\!\! \vka_3\Vert p( t)\Vert_{L^1(\Om)}\!+ \!\!\int \Big(\frac{\ga\vka_1\vka_2}2|\psi(x,t)|^{4}\!+\!C|p(x, t)|^{4}\Big) dx
\le \fr{\ga\vka_2}{2}{\cal U}(t)\!+\!C(p_0)
\eeqn
since $p(\tau)\in E \subset L^{4}(\Om)$ by the Sobolev embedding theorem in the dimension two. Moreover, \eqref{pT} gives
\be\la{fp2}
 |\langle \na \psi(t),i\na p( t)\rangle|\le \fr\ga4 \Vert\na \psi(t)\Vert^2+
 \fr1{\ga}\Vert \na p( t)\Vert^2 \le \frac{\ga}4\Vert\na \psi(t)\Vert^2 + \fr1{\ga}\,p_0.
\ee
Now \eqref{Eb2}--(\ref{fp2}) imply
\be\la{dcH}
\dot\cH (t) \le -\ga\Vert \nabla\psi(t)\Vert^2-\ga\vka_2\, {\cal U}(t)+\fr\ga4 \Vert\na \psi(t)\Vert^2+\fr{\ga\vka_2}2 {\cal U}(t)+C_1(p_0)
\le -\al_+\cH (t)+C_1(p_0),\quad t\in\R,
\ee
where $\al_+:=\fr {\ga}2\min(3,\vka_2)>0$.
 Hence the Gronwall inequality gives
\be\la{Eb4}
\cH (t)\le C(\cH(0))e^{-\al_+ t}+D(p_0),\qquad t\ge 0.
\ee
Now
 \eqref{U1}
implies \eqref{Eb5} for $t\ge 0$.
\smallskip

For $t\le 0$ the bound \eqref{Eb5} follows by similar arguments applied to the
time-reversed equation (\ref{dS}) for the function $\vp(t)=\psi(-t)$:
\be\la{dSHf2}
\dot\vp(t)=iD\cH(\vp(t))+\ga\vp(t)+ip(- t),\qquad t\in\R.
\ee
Now for $\cH_-(t):=\cH(\vp(t))$ we obtain similarly to \eqref{Eb2},
 \be\la{Eb22}
\dot\cH_- (t)=
\ga\langle \na \vp(t),\na\vp(t)\rangle+\ga\langle f(\vp(t)), \vp(t)\rangle
+\langle \na \vp(t),i\na p(-t)\rangle+\langle f(\vp(t)), ip(- t)\rangle.
\ee
In this case,  bound (\ref{Eb3}) changes to
\be\la{Eb32}
|\langle f(\vp(t)), \vp(t)\rangle|\le \frac{2\vka_3}{\vka_1}\,{\cal U}(t)+C|\Omega|,
\ee
which follows from \eqref{U1} and (\ref{U3}). Now (\ref{Eb22}) together with (\ref{fp3}), (\ref{fp2}) and (\ref{Eb32}) imply
\be\la{dcH2}
\dot\cH_- (t) \le \al_-\cH_- (t)+C,\qquad t>0,
\ee
with some $\al_->0$. Hence the Gronwall inequality gives
$
\cH_- (t)\le C(\cH(0))e^{\al_- t}+C$ for $t>0$, 
which implies \eqref{Eb5} for $t< 0$ by \eqref{U1}.
\end{proof}



\section{Well-posedness }\la{s3}
In this section, we prove the 
well-posedness of the Cauchy problem for the Schr\"odinger equation \eqref{dS} in the energy class. The key role is played by the energy equation \eqref{Eb2}.

\bt\la{pwp}
Let conditions \eqref{U1}--\eqref{U5}, \eqref{pT} hold. Then
 \smallskip\\
 {\rm i)} For any initial state $\psi(0)\in E$ equation \eqref{dS} admits a unique solution
$\psi(t)\in C(\R,E)\cap C^1(\R, E^*)$, and the bounds \eqref{Eb5} hold.
\smallskip\\
 {\rm ii)} The energy equation \eqref{Eb2} holds,
 \be\la{ebal2}
\fr d{dt}\cH (\psi(t))=\langle - \De \psi(t)+f(\psi(t)),-\ga \psi(t)-i p( t)\rangle,\qquad t\in\R.
\ee
{\rm iii)} The map $S:(\psi(0),p)\mapsto \psi(\cdot)$ is continuous
\be\la{W}
 S: E\times C(\R,E)\to C(\R,E).
 \ee
 
 \et
\subsection{Weak solutions}
First, we construct `weak solutions'.
\bl\la{lws}
Let the conditions of Theorem {\rm \ref{pwp}} hold. Then, for any initial state $\psi(0)\in E$, 
equation \eqref{dS} admits a unique solution $\psi(t)\in C(\R,E_w)\cap L^\infty_\loc(\R,E)\cap W^{1,\infty}_\loc(\R, E^*)$,
 and the bounds \eqref{Eb5} hold.
\el
\begin{proof}
The existence of a weak solution $\psi(t)\in C(\R,E_w)\cap L^\infty_\loc(\R,E)\cap W^{1,\infty}_\loc(\R, E^*)$
follows by the Galerkin approximations.
We recall this construction in Appendix \ref{A} since we will use it in the proof of \eqref{ebal2}.
The uniqueness of this solution in the case $p(x,t)\equiv 0$ and $\ga=0$
is deduced in Theorem 3.6.1 of \ci{C2003} from the Trudinger inequality \ci{A1975}.
The proof of the uniqueness for $\ga> 0$ and $p(x,t)$ satisfying \eqref{pT} remains almost unchanged:
we give the required modifications in Appendix \ref{B}. The uniqueness implies that the convergence (\ref{AA})
holds actually for the entire sequence of Galerkin approximations: for any $s<1$
\be\la{AA0}
\psi_{m}(t) \toCHs \psi(t),\qquad m\to\infty,\quad T>0.
\ee
The uniform bounds (\ref{Eb5m}) hold for the Galerkin approximations $\psi_m$. Hence the bounds \eqref{Eb5}
hold also for their limit $\psi(t)$.
\epr
\subsection{Energy equation}
Let us prove the energy equation \eqref{ebal2}.
\bp\la{pen}
Let the conditions of Theorem {\rm \ref{pwp}} hold. Then the energy equation \eqref{ebal2} holds for any solution
 $\psi(t)\in C(\R,E_w) \cap L^\infty_\loc(\R,E)\cap W^{1,\infty}_\loc(\R, E^*)$ to \eqref{dS}.
 \ep
\bpr
The identity \eqref{ebal2} for smooth solutions can be rewritten as
$$
\frac 12\pa_t\Vert\na\psi(t)\Vert^2+\ga\Vert\psi(t)\Vert^2+\pa_t \cU(\psi(t))
=\langle\na \psi(t), -i\na p( t)\rangle+\langle f(\psi(t)),-\ga \psi(t)-i p( t)\rangle,\qquad t\in\R.
$$
Equivalently,
\beqn\la{ebal3}
\!\!\!\!\!\!\!\!\!\!\!\!
\frac 12[e^{2\ga t}\langle \na\psi(t), \na\psi(t)\rangle] \Big|_0^T\!\!&\!\!=\!\!&\!\!-[e^{2\ga t}\cU(\psi(t))]\Big|_0^T
+2\ga\int_0^T e^{2\ga t}\cU(\psi(t))dt
\nonumber\\\nonumber\\
&&\!\!-\int_0^Te^{2\ga t}\langle f(\psi(t)),\!\ga \psi(t)\!+\!i p( t)\rangle dt-\int_0^Te^{2\ga t}\langle\na \psi(t), i\na p( t)\rangle dt,\quad T\in\R.
\eeqn
We will deduce \eqref{ebal3} by the limit transition in the corresponding equation for the Galerkin approximations.
Namely, the Galerkin equations (\ref{dSm}) imply for any $T\in\R$
\beqn\la{ebal5}
\frac 12[e^{2\ga t}\langle \na\psi_m(t), \na\psi_m(t)\rangle] \Big|_0^T\!\!\!&\!=\!&\!\!\!-[e^{2\ga t}\cU(\psi_m(t))]\Big|_0^T
+ 2\ga\int_0^T\!e^{2\ga t}\cU(\psi_m(t)) dt
\nonumber\\\nonumber\\
&&\!\!\!-\int\limits_0^T\!e^{2\ga t}\langle f(\psi_m(t)),\!\ga \psi_m(t)\!+\!i p_m( t)\rangle dt-\int\limits_0^Te^{2\ga t}\langle\na \psi_m(t), i\na p_m( t)\rangle dt.
\eeqn
\bl\la{lrhs}
{\rm a)} The right-hand side of \eqref{ebal5} converges to the same of \eqref{ebal3}.
\\
{\rm b)} The left-hand side of \eqref{ebal3} is dominated by the limit of the left-hand side of \eqref{ebal5}.
\el
\begin{proof}
It suffices to consider the case $T>0$.
\smallskip\\
a) The convergence (\ref{AA0}) of the Galerkin approximations with any $s<1$ implies that
\be\la{AAp}
\psi_{m}(t) \toCLr\psi(t),\qquad m\to\infty.
\ee
with any $r<\infty$ by the Sobolev embedding theorem in the dimension two.
Hence the convergence of the first three terms on the right-hand side of \eqref{ebal5} to the same of \eqref{ebal3} follows.
\smallskip

The last term of \eqref{ebal5} converges to the last term of \eqref{ebal3} by the Lebesgue theorem due to
1) the uniform bounds (\ref{Eb5m}),
2) the strong convergence $p_m(t)\toE p(t)$ and
3) the weak convergence
$\psi_m(t)\toEw \psi(t)$ for all $t\ge 0$. This weak convergence follows from (\ref{AA0}) and (\ref{Eb5m}).\\
b) The second assertion of the lemma follows from the weak convergence
$\psi_m(T)\toEw \psi(T)$, and the strong convergence $\psi_m(0):=P_m\psi(0)\toE \psi(0)$.
\end{proof}
As a result, we obtain in the limit the {\it energy inequality}: the left-hand side of \eqref{ebal3} is dominated by its right-hand side for every $T\in\R$.

On the other hand, $\psi(t)\in C(\R,E_w)$ by Lemma \ref{lws}, and hence, $\psi(T)\in E$. Therefore, the same lemma implies the existence
of a~solution $\ti\psi(t)\in C(\R,E_w)\cap L^\infty_\loc(\R,E)\cap W^{1,\infty}_\loc(\R, E^*)$ starting from 
the initial state $\ti\psi(T)=\psi(T)\in E$.
The same arguments of Lemma \ref{lrhs} show that the right-hand side of \eqref{ebal3}, with $\ti\psi$ instead of $\psi$,
is dominated by its left-hand side.

At last, the uniqueness of a~solution in Lemma \ref{lws} implies that $\ti\psi(t)=\psi(t)$ for $t\in\R$.
Hence the integral energy equations \eqref{ebal3} holds for $t\in\R$.
Finally, \eqref{ebal2} holds for a.a. $t\in\R$ and in the sense of distributions.
The proposition is proved.
\epr

\subsection{Strong solutions}

Now we can prove that $\psi(t)\in C(\R,E)\cap C^1(\R,E^*)$. Namely, $\psi(t)\in C(\R,E_w)\cap L^\infty_\loc(\R,E)$ by Lemma \ref{lws}.
On the other hand, the energy equation \eqref{ebal3} holds for all $T\in\R$, and
the right-hand side of \eqref{ebal3} is a continuous function of $T\in\R$ by the same arguments as in the proof of Lemma \ref{lrhs}.
Hence \eqref{ebal3} implies that the norm $\Vert\psi(t)\Vert_E$ is also a~continuous function of $t\in\R$. Therefore, $\psi(t)\in C(\R, E)$, and
the equation \eqref{dS} together with Lemma \ref{lnl},~{\rm ii)} and condition \eqref{pT} imply that $\dot\psi(t)\in C(\R, E^*)$.

\subsection{Continuous dependence}
 The continuity of the map (\ref{W}) follows by similar arguments. Namely, let $\psi(0), \psi^n(0)\in E$, $p,p^n\in C(\R, E)$ and
\be\la{wpp2}
\sup_n\Vert\psi^n(0)\Vert_E<\infty,\quad \psi^n(0) \toE \psi(0); \quad p^n \toCE p,\qquad n\to\infty.
\ee
Consider the corresponding unique solutions $\psi^n(t)\in C(\R,E)\cap C^1(\R,E^*)$ to \eqref{dS}. Similarly to (\ref{AA0}), we obtain for any $s<1$ and any $T>0$
\be\la{AA1}
\psi^n(t) \toCHs\psi(t),\qquad n\to\infty.
\ee
Hence the Banach theorem on weak compactness implies that
$$
\psi^n(t)\toEw\psi(t),\qquad n\to\infty,\qquad t\in\R.
$$
On the other hand, the energy equation \eqref{ebal3} implies the norm-convergence
$$
\Vert\psi^n(t)\Vert_E \toC\Vert\psi(t)\Vert_E,\qquad n\to\infty
$$
 by the same arguments as in the proof of Lemma \ref{lrhs}, using the uniform bounds \eqref{Eb5} and the uniform convergence (\ref{AA1}). Hence
 $$
\psi^n(t)\toCE\psi(t),\qquad n\to\infty.
$$

Now Theorem \ref{pwp} is proved.
\subsection{Weak continuity of the process}
In conclusion, we prove a technical lemma which will be important in the next section.
\bl\la{lw}
Let conditions \eqref{U1}--\eqref{U4} and \eqref{pT} hold, and let $1/3\le s<1$.
Then, for any $T>0$, the map $S(T,0):(\psi(0),p)\mapsto \psi(T)$ is continuous
from $H^s\times C([0,T],E)$ to $H^s$.
\el
\begin{proof}
Let $\psi(0), \psi^n(0)\in E$, $p,p^n\in C([0,T],E)$ and
\be\la{wpp}
\sup_n\Vert\psi^n(0)\Vert_E<\infty,\quad \psi^n(0) \toHs \psi(0); \quad p^n \toCTE p,\qquad n\to\infty.
\ee
Consider the corresponding unique solutions $\psi^n(t)\in C(\R,E)\cap C^1(\R,E^*)$ to
\be\la{dSn}
i\dot\psi^n(x,t)=-\De\psi^n(x,t)-i\ga\psi^n(x,t)+f(\psi^n(x,t))+p^n(x, t),\qquad x\in\Om,\qquad t\in\R.
\ee
These trajectories are equicontinuous on $[0,T]$ with values in $H^s$. This follows from the interpolation inequality (\ref{intin})
since $\sup_n \sup_{t\in [0,T]}\Vert\psi^n(t)\Vert_{H^1}<\infty$ by \eqref{Eb5},
while $\sup_n \sup_{t\in [0,T]}\Vert\dot \psi^n(t)\Vert_{E^*}<\infty$ by equations (\ref{dSn}) together with Lemma \ref{lnl},~ii) and \eqref{pT}.

Now the Arzel\`a--Ascoli theorem implies that there exists a subsequence $\psi^{n'}$ converging in $C([0,T],H^s)$:
\be\la{AA3}
\psi^{n'} \toCnHs \ti\psi,\qquad n'\to\infty.
\ee
The limit function $\ti\psi(t)\in C([0,T],H^s)$ satisfies the Schr\"odinger equation \eqref{dS}. This follows by
the limit transition in the equation (\ref{dSn}) taking into account Lemma \ref{lnl},~i)
 and the uniform convergence of $p^n$ in (\ref{wpp}). Moreover, $\ti\psi(0)=\psi(0)$ by the convergence of $\psi^n(0)$ in (\ref{wpp}).
Hence the uniqueness arguments as above imply that $\ti\psi(t)\equiv\psi(t)$ for all $t\in [0,T]$.
Therefore, the continuity of $S(T,0)$ follows from (\ref{AA3}) since the limit does not depend on the subsequence $n'$.
\end{proof}
\section{Global attractor}\la{s4}
In this section, we prove the convergence
in the $H^1$-norm
of all finite energy solutions to a compact attractor.
Proposition \ref{pwp},~i) allows us to define the continuous process $S_p(t,\tau)$ in $E$ acting by
\be\la{Wtt}
S_p(t,\tau):\psi(\tau)\mapsto \psi(t),
\ee
where $\psi(t)\in C(\R,E)\cap C^1(\R, E^*)$ is a solution to \eqref{dS}.
The bounds \eqref{Eb5} imply that the process admits a~bounded uniform absorbing set in $H^1$.
Namely, let $B_R$ denote the ball $\{\psi\in E:\Vert\psi\Vert_E\le R\}$.
Denote $\hat B:=B_{\sqrt{D+1}}$, where $D$ is the constant from \eqref{Eb5}.
Then \eqref{Eb5} implies that
$$
S_p(t,\tau)B_R\subset \hat B,\qquad t-\tau> \fr{\log A(R)}\alpha.
$$
Now we prove the existence of a uniform compact attractor in $H^1$.
Let us recall the definition of the uniform attractor for the nonautonomous
equations, see, e.g., \ci{Har1991}, \cite[Definition A2.3]{CV1992}, \cite[Definition (7)]{CV1993}.
\bd\la{dA} \rm 
A set $\cA\subset E$ is called a~\textit{uniform attractor} of the process $S_p(t,\tau)$
if it is a uniformly attracting set and if any other closed attracting set contains $\cA$.
\ed
The term `uniformly attracting' means that, for any bounded subset $B\subset E$,
\be\la{attr}
\dist_E (S_p(t,\tau)B,\cA)\to 0,\quad t-\tau\to\infty.
\ee
\bt \la{tA}
Let the assumptions of Theorem {\rm \ref{pwp}} hold,
$\ga>0$, and let the pumping $p(t)$ be an almost periodic function with values in $E$.
Then the process $S_p(t,\tau)$ admits a compact uniform attractor in $E$.
\et
For the proof we develop the strategy \ci{B2004} for a nonautonomous equation.
Namely, we define the functionals (cf. (2.3) and (2.4) of \ci{L1995} and (2.12), (2.25) of \ci{W1995})
 \be\la{func}
 \Phi(\psi,t):=\cH(\psi)+\langle p( t),\psi \rangle,\qquad
 \Psi(\psi,t):=\cU(\psi)-\fr12\langle f(\psi),\psi \rangle+\fr12\langle p( t),\psi \rangle,\quad \psi\in E,\quad t\ge 0.
 \ee
 \bl \la{lwp}
 For any $s\ge 1/2$ and any $t\ge 0$, the functionals $\cU(\cdot)$ and $\Psi(\cdot,t)$ are continuous on $H^s$ and bounded on bounded subsets.
 \el
 \begin{proof}
 Using (\ref{U3}) and \eqref{U4}, we obtain
 $$
 |\cU(\psi_1)-\cU(\psi_2)|\le C\int_\Om (1+|\psi_1(x)|^3+|\psi_2(x)|^3)|\psi_1(x)-\psi_2(x)|dx,
$$
\beqn\nonumber
&&|\langle f(\psi_1),\psi_1\rangle-\langle f(\psi_2),\psi_2 \rangle|\le |\langle f(\psi_1)-f(\psi_2),\psi_1\rangle|
+|\langle f(\psi_2),\psi_1-\psi_2 \rangle|\\
\nonumber
&&\le C\Big(\int_\Om (1+|\psi_1(x)|^2+|\psi_2(x)|^2)|\psi_1(x)-\psi_2(x)||\psi_1|dx+ \int_\Om (1+|\psi_2(x)|^3)|\psi_1(x)-\psi_2(x)|dx\Big)\\
\nonumber
 &&\qquad\qquad\qquad\qquad\qquad\qquad\qquad\qquad\qquad \le C_1\int_\Om (1+|\psi_1(x)|^3+|\psi_2(x)|^3)|\psi_1(x)-\psi_2(x)|dx
\eeqn
Applying the H\"older inequality and the Sobolev embedding theorem, we obtain for any $s\ge 1/2$
\begin{align}\nonumber
 |\cU(\psi_1)-\cU(\psi_2)|+|\langle f(\psi_1),\psi_1\rangle-\langle f(\psi_2),\psi_2 \rangle|
 &\le C(1+\Vert\psi_1\Vert_{L^4}^3+\Vert\psi_3\Vert_{L^4}^3)\Vert\psi_1-\psi_2\Vert_{L^4}\\
 \nonumber
 &\le C(1+\Vert\psi_1\Vert_{H^s}^3+\Vert\psi_3\Vert_{H^s}^3)\Vert\psi_1-\psi_2\Vert_{H^s}.
\end{align}
Similarly,
$$
| \langle p(t),\psi_1-\psi_2 \rangle|\le \Vert p(t)\Vert \Vert \psi_1-\psi_2\Vert\le Cp_0\Vert\psi_1-\psi_2\Vert_{H^s},
$$
for any $s\ge 0$.
\end{proof}

 Let $\psi(t)\in C(\R,E)\cap C^1(\R,E^*)$ be a solution to \eqref{dS}.
Now the energy equation \eqref{ebal2} implies that
\be\la{func2}
\fr d{dt}\Phi(\psi(t),t)=-2\ga \Phi(\psi(t),t)+2\ga \Psi(\psi(t),t)+\langle \dot p(t),\psi(t) \rangle,\qquad t>0,
 \ee
 where the last term is well-defined in the sense of distributions:
for any test function $\vp\in C^\infty(\ov{\R_+}) $
 \be\la{last}
 \int_0^t \langle \dot p(s),\psi(s) \rangle\vp(s)ds:=[\langle p(s),\psi(s)\rangle\vp(t)]
 \Big|_0^t-\int_0^t[\langle p(s),\dot\psi(s)\rangle\vp(s)+\langle p(s),\psi(s)\rangle\dot\vp(s)]ds.
 \ee
 To prove (\ref{func2}), we
 differentiate
 using
 \eqref{dS} and \eqref{ebal2},
\begin{align}\nonumber
\fr d{dt}\Phi(\psi(t),t)-\langle \dot p(t),\psi(t) \rangle&= 
\dot\cH (\psi(t))+\langle \dot\psi(t), p(t) \rangle\\
\nonumber
&=-\ga\langle \na \psi(t),\na\psi(t)\rangle-\ga\langle f(\psi(t)), \psi(t)\rangle+\langle i\na \psi(t),\na p(t)\rangle+\langle if(\psi(t)), p(t)\rangle\\
\nonumber
&\quad\,+\langle i\De\psi(t)-\ga\psi(t)-if(\psi(t))-ip(t), p(t) \rangle\\
\nonumber
&=-2\ga \cH(\psi(t))+2\ga \cU(\psi(t)) -\ga\langle f(\psi(t)), \psi(t)\rangle-\ga\langle \psi(t), p(t) \rangle\\
\nonumber
&= -2\ga \Phi(\psi(t),t)+2\ga \Psi(\psi(t),t).
\end{align}
 Integrating (\ref{func2}), we get
 \be\la{func22}
\Phi(\psi(t),t)=e^{-2\ga t} \Phi(\psi(0),0)+\int_0^t e^{-2\ga(t-\tau)}\big[ 2\ga\Psi(\psi(\tau),\tau)+\langle \dot p(\tau),\psi(\tau) \rangle\big]d\tau,\qquad t\ge 0.
 \ee

\subsection{Strong convergence}
The following proposition is the key step in the proof of Theorem \ref{tA}.
\bp\la{pA}
For any sequences $ t_k,\tau_k\to\infty$ with $t_k-\tau_k\to\infty$ and $\phi_k\in \hat B$, there is a subsequence $k^*$ and $\phi\in E$ such that
\be\la{phi2}
S_p(t_{k^*},\tau_{k^*})\phi_{k^*}\toE \phi,\qquad k^*\to\infty.
\ee
\ep
\begin{proof}
The sequence $S_p(t_k,\tau_k)\phi_k$ is bounded in $E$ by \eqref{Eb5}.
Let us fix an arbitrary $s\in[1/2,1)$.
The inclusion $E\subset H^s$ is compact
 by the Sobolev embedding theorem.
Hence there is a subsequence strongly converging in $H^s$:
\be\la{phi3}
S_p(t_{k'},\tau_{k'})\phi_{k'}\toHs \phi,\qquad k'\to\infty.
\ee
Now to prove (\ref{phi2}) it suffices to check that for a subsequence $\{k^*\}
\subset \{k'\}$
\be\la{phi}
\limsup_{k^*\to\infty}\Vert S_p(t_{k^*},\tau_{k*})\phi_{k^*}\Vert_E\le \Vert\phi\Vert_E.
\ee
Recall that $t_k-\tau_k\to\infty$.
Hence similarly to (\ref{phi3}), for any $T>0$, there exists a subsequence
 $\{k''\}\subset\{k'\}$ and an element $\phi(-T)\in E$ such that
\be\la{phi3T}
S_p(t_{k''}-T,\tau_{k''})\phi_{k''}\toHs \phi(-T),\qquad k''\to\infty.
\ee
We set  $p_{k''}(t):=p(t_{k''}-T+t)$ and denote
\be\la{psi''}
\psi_{k''}(t):=S_p(t_{k''}-T+t,\tau_{k''})\phi_{k''}=S_{p_{k''}}(t,0)S_p(t_{k''}-T,\tau_{k''})\phi_{k''},\qquad t\in [0,T].
\ee
In particular, for $t=T$
 the convergence (\ref{phi3})
implies that
\be\la{phi50}
\psi_{k''}(T):=S_p(t_{k''},\tau_{k''})\phi_{k''}\toHs \phi,\qquad k''\to\infty.
\ee

Now we are going to apply Lemma \ref{lw} iv). At this moment we need the almost-periodicity (\ref{alper}) of the pumping $p(t)$.
 It implies that for a subsequence $\{k^*\}\subset\{k''\}$
\be\la{alm}
p_{k^*} \toCbE p^*,\qquad k^*\to\infty.
\ee
Hence (\ref{phi3T}), (\ref{psi''}) and the continuity
of the map (\ref{W})
 imply the convergence
for $t\in[0,T]$,
\be\la{phi4}
\psi_{k^*}(t)=S_{p_{k^*}} (t,0)S_p(t_{k^*}-T,\tau_{k*})\phi_{k*}\toHs\psi^*(t):=S_{p^*}(t,0)\phi(-T),\quad k^*\to\infty.
\ee
In particular, for $t=T$
the convergence (\ref{phi50})
gives that
\be\la{phi5}
\psi^*(T):=S_{p^*}(T,0)\phi(-T)=\phi.
\ee
Now we apply the integral identity
 (\ref{func22}) to solutions $\psi_{k^*}(t)$ and $\psi^*(t)$
of the equation \eqref{dS} with the pumping $p_{k^*}(t)$ and $p^*(t)$ respectively.
We obtain
\begin{align}\la{func3}
\Phi(\psi_{k^*}(T),T)&=e^{-2\ga T}\Phi(\psi_{k^*}(0),0)+\int_0^T e^{2\ga (T-\tau)}[2\ga\Psi(\psi_{k^*}(\tau),\tau)+
\langle \dot p_{k^*}(\tau),\psi_{k^*}(\tau) \rangle]d\tau
\\
\nonumber\\
\la{func4}
\Phi(\psi^*(T),T)&=e^{-2\ga T}\Phi(\psi^*(0),0)+\int_0^T e^{2\ga (T-\tau)}[2\ga\Psi(\psi^*(\tau),\tau)+
\langle \dot p^*(\tau),\psi^*(\tau) \rangle]d\tau,
\end{align}
where the integrals with
 $\dot p_{k^*}(\tau)$ and $\dot p(\tau)$ are definite according to (\ref{last}).
Making $k^*\to\infty$ in (\ref{func3}) we obtain
\be\la{func5}
\limsup_{k^*\to\infty} \Phi(\psi_{k^*}(T),T)\le C e^{-2\ga T}+\int_0^T e^{2\ga (T-\tau)}[2\ga\Psi(\psi^*(\tau),\tau)+
\langle \dot p^*(\tau),\psi^*(\tau) \rangle]d\tau
\ee
by the following arguments:
\smallskip

i) $\sup_{k^*}\Vert\psi_{k^*}(0)\Vert_E<\infty$, while the functionals $\Phi(\cdot,0)$ are
bouded on bounded subsets of $E$,
\smallskip

ii) $\psi_{k^*}(\tau)\toHs \psi^*(\tau)$ by (\ref{phi4}), while the functionals $\Psi(\cdot,\tau)$ are continuous
on $H^s$ by Lemma \ref{lwp} since $s\in[1/2,1]$,
 and
 \smallskip

 iii) (\ref{alm}) implies the convergence of the integrals with
 $\dot p_{k^*}(\tau)$ by their definition (\ref{last}).
 \smallskip\\
 Now eliminating the integral term from (\ref{func4}) and (\ref{func5}), we obtain
\be\la{func51}
\limsup_{k^*\to\infty} \Phi(\psi_{k^*}(T),T)\le C e^{-2\ga T}+\Phi(\psi^*(T),T)-e^{-2\ga T}\Phi(\psi^*(0),0).
\ee
By definitions (\ref{Ham}) and (\ref{func}),
\begin{align}
\Phi(\psi_{k^*}(T),T)&=\fr12 \Vert \psi_{k^*}(T) \Vert_E
+\cU(\psi_{k^*}(T))+\langle p_{k^*}(T),\psi_{k^*}(T) \rangle,
\nonumber\\\nonumber\\
\Phi(\psi^*(T),T)&=\fr12 \Vert \psi^*(T) \Vert_E
\,\,+\,\cU(\psi^*(T))\,+\,\langle p^*(T),\psi^*(T) \rangle,
\nonumber
\end{align}
Now Lemma \ref{lwp} together with (\ref{phi3}), (\ref{phi50}) and (\ref{phi5}) imply that
$$
\lim_{k^*\to\infty} \cU(\psi_{k^*}(T))=\cU(\phi)=\cU(\psi^*(T)),
\quad \lim_{k^*\to\infty} \langle p_{k^*}(T),\psi_{k^*}(T) \rangle= \langle p^*(T),\phi\rangle=\langle p^*(T),\psi^*(T)\rangle.
$$
Hence (\ref{func51}) becomes
$$
\fr12\limsup_{k^*\to\infty} \Vert \psi_{k^*}(T) \Vert_E
\le C e^{-2\ga T}+
\fr12\Vert\phi\Vert_E
-e^{-2\ga T}\Phi(\psi^*(0),0).
$$
According to the definition (\ref{phi50}), we can replace $\psi_{k^*}(T)$ by $S_p(t_{k^*},\tau_{k*})\phi_{k^*}$. Now making 
 $T\to\infty$, we get (\ref{phi}) since $\ga>0$.
 \end{proof}
\subsection{Compactness}
Let us denote by $\cA$ the set of all points $\phi$ from (\ref{phi2}). This set is obviously closed.
Now Theorem \ref{tA} will follow from the next lemma.
\bl
{\rm i)} The set $\cA$ is compact in $E$,
and
{\rm ii)} the set $\cA$ is uniformly attracting in $E$.
\el
\begin{proof}
i) Let  us consider a sequence $\phi_n\in\cA$. Then, for each $n$, 
\be\la{phi22}
S_p(t_{nk},\tau_{nk})\phi_{nk}\toE \phi_n,\qquad t_{nk}-\tau_{nk}\to\infty
\ee
as $k\to\infty$,
where $\phi_{nk}\in \hat B$. Hence there exists a sequence $k(n)$ such  that
\be\la{phi23}
\Vert S_p(t_{nk(n)},\tau_{nk(n)})\phi_{nk(n)}- \phi_n\Vert_E\to 0,
\qquad t_{nk(n)}-\tau_{nk(n)}\to\infty
\ee
as $n\to\infty$.
However, Proposition \ref{pA} implies that for a subsequence $n'$
\be\la{phi24}
S_p(t_{n'k(n')},\tau_{n'k(n')})\phi_{n'k(n')}\toE \ov\phi,\qquad n'\to\infty.
\ee
Hence $\ov\phi\in\cA$, and (\ref{phi23}) implies that $\phi_{n'}\toE \ov\phi$.
Therefore, the first assertion of the lemma is proved.
\smallskip\\
ii) Let us assume the contrary. Then there exists  a~sequence 
$S(t_k,\tau_k)\phi_k$ such that
\be\la{seq}
\dist_E(S(t_k,\tau_k)\phi_k,\cA)\ge\ve>0,\qquad t_k-\tau_k\to\infty.
\ee
On the other hand, Proposition \ref{pA} implies that for a subsequence
$k'$,
\be\la{seq2}
S(t_{k'},\tau_{k'})\phi_{k'}\toE \phi\in\cA,\qquad k'\to\infty.
\ee
This contradiction proves the second assertion of the lemma.
\end{proof}

\appendix

\setcounter{section}{0}
\setcounter{equation}{0}
\protect\renewcommand{\thesection}{\Alph{section}}
\protect\renewcommand{\theequation}{\thesection.\arabic{equation}}
\protect\renewcommand{\thesubsection}{\thesection.\arabic{subsection}}
\protect\renewcommand{\thetheorem}{\Alph{section}.\arabic{theorem}}

\section{Galerkin approximations}\la{A}

Denote by $E^m$, $m=1,2,...$, the linear span of the first $m$ eigenfunctions of the
 Laplacian $\De$ on the region $\Om$ with the Dirichlet boundary conditions (\ref{Dbc}).
 We define the Galerkin approximations $\psi_m(t)$ as solutions to finite-dimensional dissipative equations in $E^m$
 \be\la{dSm}
i\dot\psi_m(t)=-\De\psi_m(t))-i\ga \psi_m(t)+P_m f(\psi_m(t))+p_m( t),\quad t\in\R; \quad \psi_m(0)=P_m\psi(0).
\ee
Here $P_m$ stands for the orthogonal projection of $L^2$ onto $E^m$, and $p_m( t):=P_m p( t)$.
Applying the calculations from the proof of Lemma \ref{l1} to the equations (\ref{dSm}),
we get the uniform estimates of type \eqref{Eb5} for $\psi_m$ with the same constants:
\be\la{Eb5m}
\Vert\psi_m(t) \Vert_E^2 \le C(\Vert\psi(0)\Vert_E)e^{-\al_\pm t}+D,\qquad t\in\R.
\ee
Hence, the Galerkin approximations $\psi_m(t)$ exist globally in time.
\smallskip

Let us show that these uniform estimates imply the existence of a limiting function $\psi(t)$ of the
Galerkin approximations.
For this purpose we will use the known continuity property of the nonlinear term:
\bl\la{lnl}
Let nonlinear function $f(\psi)$ satisfies condition \eqref{U4}. Then\\
{\rm i)}
For any $r\ge 1$ and $s\ge 1-\ds\frac{2}{3r}$,
the nonlinearity $N:\psi(\cdot)\mapsto f(\psi(\cdot))$ is the continuous map $H^s\to L^r$.
\smallskip\\
{\rm ii)} $N:\psi(\cdot)\mapsto f(\psi(\cdot))$ is the continuous map $E\to E^*$.
\smallskip\\
{\rm iii)} $N$ is the continuous map $E_1\to L^2$.

\el
\begin{proof}
i) Condition \eqref{U4} implies that
$$
|f(\psi_1(x))-f(\psi_2(x))|\le C(1+|\psi_1(x)|^{2}+|\psi_2(x)|^{2})|\psi_1(x)-\psi_2(x)|.
$$
Hence, the H\"older inequality and the Sobolev embedding theorem give, for $r\ge 1$, 
\begin{align}\nonumber
\Vert f(\psi_1)-f(\psi_2)\Vert_{L^r} &\le C\Big(1+\Vert\psi_1\Vert_{L^{3r}}^{2}+\Vert\psi_1\Vert_{L^{3r}}^{2})\Vert\psi_1-\psi_2\Vert_{L^{3r}}
\\
\la{Hin}
&\le C(s)(1+\Vert\psi_1\Vert_{H^s}^{2}+\Vert\psi_2\Vert_{H^s}^{2}) \Vert\psi_1-\psi_2\Vert_{H^s}.
\end{align}

ii) In particular, (\ref{Hin}) implies that $N$ is continuous $H^{2/3}\to L^2$.
It remains to note that $E\subset H^{2/3}$, while $L^2\subset E^*$.
\smallskip\\
iii) The map
$N:E_1\to C(\ov\Om)$ is continuous by the Sobolev embedding theorem.
\end{proof}

Further, the Galerkin approximations are equicontinuous in $H^s$ with any $s<1$ by the Dubinsky `Theorem on Three Spaces' ( \ci[Theorem 5.1]{L1969}).
Namely, this equicontinuity follows from the 
interpolation inequality: for any $\de>0$
\be\la{intin}
\Vert\psi_m(t_1)-\psi_m(t_2)\Vert_{H^s}\le \de \Vert\psi_m(t_1)-\psi_m(t_2)\Vert_{H^1}+C_\de
\Vert\psi_m(t_1)-\psi_m(t_2)\Vert_{E^*},\qquad t_1,\,t_2\in\R.
\ee
Here the first term on the right is small for sufficiently small $\de>0$, 
since $\sup_m \sup_{t\ge 0}\Vert\psi_m(t)\Vert_E<\infty$ by (\ref{Eb5m}),
while the second term is small for $|t_1-t_2|\ll 1$, 
since $\sup_m \sup_{t\ge 0}\Vert\dot \psi_m(t)\Vert_{E^*}<\infty$ by the Galerkin equations (\ref{dSm}) together with Lemma \ref{lnl},~ii) and \eqref{pT}.

Now the Arzel\`a--Ascoli theorem implies that there exists a subsequence $\psi_{m'}$
converging in $C(\R,H^s)$:
\be\la{AA}
\psi_{m'}(t) \toCHs\psi(t),\qquad m'\to\infty.
\ee
The limit function $\psi(t)\in C(\R,H^s)$ satisfies the
 estimates \eqref{Eb5} by the uniform bounds (\ref{Eb5m}),
 and hence, $\psi(t)\in C(\R,E_w)\cap L^\infty_\loc (\R,E)$.

The limit function $\psi(t)$ satisfies the Schr\"odinger equation \eqref{dS} in the sense of distributions.
This follows from the convergence (\ref{AA}) with $s<1$ close to $1$ by limit transition in the Galerkin equations (\ref{dSm})
taking into account Lemma \ref{lnl}.
Finally, the convergence $\psi_m(0)\to \psi(0)$ as $m\to\infty$ implies that $\psi(0)=\psi(0)$.

At last, the equation \eqref{dS} implies that $\psi(t)\in W^{1,\infty}_\loc(\R,E^*)$ by Lemma \ref{lnl},~ii).

\section{Uniqueness}\la{B}

Here we reduce the proof of the uniqueness of solution to \eqref{dS}
with $\ga>0$ and fixed $\psi(0)\in E$
to the case $\ga=0$, which was considered in Section 3.6 of \ci{C2003}.
Namely, let $\psi_1,\psi_2\in L^\infty_\loc(\R,E)\cap W^{1,\infty}_\loc(R, E^*)$
be two solutions to \eqref{dS}
with the same initial state $\psi_1(0)=\psi_2(0)$.
Then the difference $z=\psi_1-\psi_2$ satisfies the equation
$$
i\dot z(x,t)=-\De z(x,t)-i\ga z(x,t)+f(\psi_1(x,t))-f(\psi_2(x,t)),\qquad x\in\Om,\quad t\in\R.
$$
This equation implies
\be\la{dSa2}
\fr {d}{dt}\Vert z(t)\Vert_{L^2}^2+2\ga \Vert z(t)\Vert_{L^2}^2 =
2\langle f(\psi_1(t))-f(\psi_2(t)),z(t)\rangle,\quad t\in\R.
\ee
Equivalently, we have the equation
\be\la{dSa3}
\fr {d}{dt}[e^{2\ga t}\Vert z(t)\Vert_{L^2}^2] =2e^{2\ga t} \langle f(\psi_1(t))-f(\psi_2(t)),z(t)\rangle,\quad t\in\R.
\ee
Finally, denoting $h(t)=|\psi_1(t)|+|\psi_2(t)|$,
$w(t)=e^{\ga t}z(t)$ and using \eqref{U4}, we obtain the inequality
\be\la{dSa4}
|\fr {d}{dt}\Vert w(t)\Vert_{L^2}^2|\le
C\int_\Om (1+h^2(t)) | w(x,t)|^2dx,\quad t\in\R,
\ee
which
coincides with the corresponding inequality in Section 3.6 of \ci{C2003}.
Using further the Trudinger inequality as in \ci{C2003}, we obtain  $w(t)\equiv 0$.


\end{document}